\documentclass{amsart}
\usepackage{amsmath,amssymb}
\newtheorem{theorem}{Theorem}
\newtheorem{lemma}[theorem]{Lemma}
\newtheorem{proposition}[theorem]{Proposition}
\newtheorem{corollary}[theorem]{Corollary}

\newtheorem{example}[theorem]{Example}

\begin{document}
\title{Localization of injective modules over valuation rings}
\author{Fran\c cois Couchot}
\address{Laboratoire de Math\'ematiques Nicolas Oresme, CNRS UMR
  6139,
D\'epartement de math\'ematiques et m\'ecanique,
14032 Caen cedex, France}
\email{couchot@math.unicaen.fr} 



\begin{abstract} It is proved that $E_J$ is injective if $E$ is an injective module over a valuation ring $R$,
for each prime ideal $J\ne Z$. Moreover, if $E$ or $Z$ is flat, then $E_Z$ is injective too. It follows that localizations of
injective modules over h-local Pr\"ufer domains are injective too.
\end{abstract}
\maketitle

If $S$ is a multiplicative subset of a noetherian ring $R$, it is well known that $S^{-1}E$ is injective for each injective $R$-module $E$.
The following example shows that this result is not generally true if $R$ is not noetherian.
\begin{example} \label{E:VN} \textnormal{Let $K$ be a field and $I$ an infinite set. We put $R=K^I$, $J=K^{(I)}$ and $S=\{1-r\mid r\in J\}$. Then $R/J\cong S^{-1}R$,
$R$ is an injective module, but $R/J$ is not injective by \cite[Theorem]{Osof}.}
\end{example}

However, we shall see that, for some classes of non-noetherian rings, localizations of injective modules are injective too. For instance:

\begin{proposition} \label{P:hered}
Let $R$ be a hereditary ring. For each multiplicative subset $S$ of $R$ and for every injective $R$-module $E$, $S^{-1}E$ is injective.
\end{proposition}
There exist non-noetherian hereditary rings.

\textbf{Proof.} Let $F$ be the kernel of the natural map: $E\rightarrow S^{-1}E$. Then $E/F$ is injective and $S$-torsion-free.
Let $s\in S$. We have $(0:s)=Re$, where $e$ is an idempotent of $R$. It is easy to check that $s+e$ is a non-zerodivisor. So, if $x\in E$,
there exists $y\in E$ such that $x=(s+e)y$. Clearly $eE\subseteq F$. Hence $x+F=s(y+F)$. Therefore the multiplication by $s$ in $E/F$ is
bijective, whence $E/F\cong S^{-1}E$. \qed

\bigskip
In Proposition~\ref{P:hered} and Example~\ref{E:VN}, $R$ is a coherent ring. By \cite[Proposition 1.2]{Cou1} $S^{-1}E$ is fp-injective if $E$ is a fp-injective module over a coherent
ring $R$, but the coherence hypothesis can't be omitted: see \cite[Example p.344]{Cou1}.

\bigskip
The aim of this paper is to study localizations of injective modules and fp-injective modules over a valuation ring $R$.
Let $Z$ be the subset of its zerodivisors. Then $Z$ is a prime ideal. We will show the
following theorem:
\begin{theorem} \label{T:main} Let $R$ be a valuation ring, denote by $Z$ the set of zero divisors of $R$ and
let $E$ be an injective (respectively fp-injective) module. Then:
\begin{enumerate}
\item For each prime ideal $J\ne Z$, $E_J$ is injective (respectively fp-injective).
\item $E_Z$ is injective (respectively fp-injective) if and only if $E$ or $Z$ is flat.
\end{enumerate}
\end{theorem}

In this paper all rings are associative and commutative with unity and
all mo\-dules are unital. We say that an $R$-module
$E$ is \textbf{divisible} if, for every $r\in
R$ and $x\in E,$ $(0:r)\subseteq (0:x)$ implies that $x\in
rE$, and that $E$ is \textbf{fp-injective}(or
\textbf{absolutely pure}) if $\mathrm{Ext}_R^1(F,E)=0,$ for every finitely
  presented $R$-module $F.$ A ring $R$ is called \textbf{self
    fp-injective} if it is fp-injective as $R$-module. An exact sequence \ $0 \rightarrow F \rightarrow E \rightarrow G \rightarrow 0$ \ is \textbf{pure}
if it remains exact when tensoring it with any $R$-module. In this case
we say that \ $F$ \ is a \textbf{pure} submodule of $E$. Recall that a
module $E$ is fp-injective if and only if it is a pure submodule
of every overmodule. A module is said to be
\textbf{uniserial} if its  submodules are linearly ordered by inclusion
and a ring $R$ is a \textbf{valuation ring} if it is uniserial as
$R$-module. Recall that every finitely presented module over a
valuation ring is a finite direct sum of cyclic modules
\cite[Theorem 1]{War1}. Consequently a module $E$ over a valuation ring
$R$ is fp-injective if and only if it is divisible. 

An $R$-module $F$ is \textbf{pure-injective} if for every pure exact
sequence
\[0\rightarrow N\rightarrow M\rightarrow L\rightarrow 0\]
 of $R$-modules, the following sequence
 \[0\rightarrow\mathrm{Hom}_R(L,F)
\rightarrow\mathrm{Hom}_R(M,F)\rightarrow\mathrm{Hom}_R(N,F)\rightarrow
0\] is exact. Then a module is injective if and only if it is pure-injective and fp-injective.
A ring $R$ is said to be an \textbf{IF-ring} if every injective module is flat. By \cite[Theorem 2]{Col} $R$ is an IF-ring if and only if
$R$ is coherent and self fp-injective.

\bigskip
In the sequel $R$ is a valuation ring whose maximal ideal is $P$ and $Z$ is its subset of zerodivisors.
Some preliminary results are needed to show Theorem~\ref{T:main}.
\begin{proposition} \label{P:prinj}
Let $R$ be a valuation ring, let $E$ be an injective module and $r\in P$. Then $E/rE$ is injective over $R/rR$.
\end{proposition}
\textbf{Proof.} Let $J$ be an ideal of $R$ such that $Rr\subset J$ and
$g:J/Rr\rightarrow E/rE$ be a nonzero homomorphism. For
each $x\in E$ we denote by $\bar{x}$ the image of $x$ in
$E/rE$. Let $a\in J\setminus Rr$ such that
$\bar{y}=g(\bar{a})\ne 0$. Then $(Rr:a)\subseteq (rE:y)$. Let
$t\in R$ such that $r=at$. Thus $ty=rz$ for some $z\in E$. It
follows that $t(y-az)=0$. So, since $at=r\ne 0$, we have $(0:a)\subset
Rt\subseteq (0:y-az)$. The injectivity of $E$ implies that there
exists $x\in E$ such that $y=a(x+z)$. We put $x_a=x+z$. If $b\in
J\setminus Ra$ then $a(x_a-x_b)\in rE$. Hence $x_b\in
x_a+(rE:_{E}a)$. Since $E$ is
pure-injective, by \cite[Theorem 4]{War} there exists $x\in\cap_{a\in
  J}x_a+(rE:_{E}a)$. It follows that
$g(\bar{a})=a\bar{x}$ for each $a\in J$. \qed

\begin{lemma} \label{L:unis} Let $R$ be a valuation ring, let $U$ be a module and $F$ a
  flat module. Then, for each $r,s\in
  R$, $F\otimes_R(sU:_Ur)\cong (F\otimes_RsU:_{F\otimes_RU}r)$.  
\end{lemma}
\textbf{Proof.} We put $E=F\otimes_RU$. Let $\phi$ be the composition of the
  multiplication by $r$ in $U$ with the natural map
$U\rightarrow U/sU$. Then $(sU:_Ur)=\mathrm{ker}(\phi)$. It follows that
  $F\otimes_R(sU:_Ur)$ is isomorphic to
  $\mathrm{ker}(\mathbf{1}_{F}\otimes\phi)$ since $F$ is flat. 
We easily
  check that $\mathbf{1}_{F}\otimes\phi$ is the composition
  of the multiplication by $r$ in $E$ with the natural map
  $E\rightarrow E/sE$. It follows that
$F\otimes_R(sU:_Ur)\cong (sE:_Er)$. \qed

\begin{proposition} \label{P:compact} Let $R$ be a valuation ring. Then every pure-injective $R$-module $F$
  satisfies the following property: if $(x_i)_{i\in I}$ is a family of
  elements of $F$ and $(A_i)_{i\in I}$ a family of ideals of
  $R$ such that the family $\mathcal{F}=(x_i+A_iF)_{i\in I}$
  has the finite intersection property, then $\mathcal{F}$ has a
  non-empty intersection. The converse holds if $F$ is flat.
\end{proposition}
\textbf{Proof.} 
 Let $i\in I$ such that $A_i$ is not finitely
generated. By \cite[Lemma 29]{Cou} either $A_i=Pr_i$ or
$A_i=\cap_{c\in R\setminus A_i}cR$. If, $\forall i\in I$ such that
 $A_i$ is not finitely generated, we replace  $x_i+A_iF$ by
$x_i+r_iF$ in the first case, and by the family
$(x_i+cF)_{c\in R\setminus A_i}$ in the second case, we deduce from
$\mathcal{F}$ a family $\mathcal{G}$ which has the finite intersection
property. Since $F$ is pure-injective, it follows that there exists
$x\in F$ which belongs to each element of the family
$\mathcal{G}$ by \cite[Theorem 4]{War}. We may assume that the family
 $(A_i)_{i\in I}$ has no smallest element. So, if $A_i$ is not finitely
 generated, there exists $j\in I$ such that
$A_j\subset A_i$. Let $c\in A_i\setminus PA_j$ such that
 $x_j+cF\in\mathcal{G}$. Then $x-x_j\in cF\subseteq A_iF$ and
 $x_j-x_i\in A_iF$. Hence $x-x_i\in A_iF$ for each $i\in I$. 

Conversely, if $F$ is flat then by Lemma~\ref{L:unis} we have
$(sF:_Fr)=(sR:r)F$ for each $s,r \in R$. We use \cite[Theorem 4]{War}
to conclude. \qed

\begin{proposition} \label{P:tenspur}
Let $R$ be a valuation ring and let $F$ be a flat pure-injective module. Then:
\begin{enumerate}
\item  $F\otimes_RU$ is pure-injective if $U$ is a uniserial module.
\item For each prime ideal $J$, $F_J$ is pure-injective.
\end{enumerate}
\end{proposition}
\textbf{Proof.} 

$(1)$. Let $E=F\otimes_RU$.
We use \cite[Theorem 4]{War} to prove
that $E$ is pure-injective. Let
$(x_i)_{i\in I}$ be a family of elements of $F$ such that the family
$\mathcal{F}=(x_i+N_i)_{i\in I}$ has the finite intersection
property, where $N_i=(s_iE:_Er_i)$ and $r_i,s_i\in R$, $\forall i\in I$. 

 First we assume that $U=R/A$ where
$A$ is a proper ideal of $R$. So $E\cong F/AF$. If $s_i\notin
A$ then $N_i=(s_iF:_Fr_i)/AF=(Rs_i:r_i)F/AF$. We set
$A_i=(Rs_i:r_i)$ in this case. If $s_i\in
A$ then $N_i=(AF:_Fr_i)/AF=(A:r_i)F/AF$. We put $A_i=(A:r_i)$ in this
case. For each $i\in I$, let $y_i\in F$ such that $x_i=y_i+AF$. It is obvious that the family $(y_i+A_iF)_{i\in I}$ has the finite
intersection property. By Proposition~\ref{P:compact} this family has a 
non-empty intersection. Then $\mathcal{F}$ has a 
non-empty intersection too. 

Now we assume that $U$ is not finitely generated. It is obvious that
$\mathcal{F}$ has a non-empty intersection if $x_i+N_i=E,\ \forall i\in I$. Now assume there exists $i_0\in I$
such that $x_{i_0}+N_{i_0}\ne E$. Let $I'=\{i\in I\mid N_i\subseteq
N_{i_0}\}$ and $\mathcal{F}'=(x_i+N_i)_{i\in I'}$. Then $\mathcal{F}$ and
  $\mathcal{F}'$ have the same intersection. By
Lemma~\ref{L:unis}
$N_{i_0}=F\otimes_R(s_{i_0}U:_Ur_{i_0})$. It follows that
$(s_{i_0}U:_Ur_{i_0})\subset U$ because $N_{i_0}\ne E$.  Hence  $\exists u\in U$ such that
  $x_{i_0}+N_{i_0}\subseteq F\otimes_RRu$. Then, $\forall
  i\in I'$, $x_i+N_i\subseteq F\otimes_RRu$. We have
  $F\otimes_RRu\cong F/(0:u)F$. From the first part of the proof $F/(0:u)F$ is pure-injective. So we may replace $R$ with
  $R/(0:u)$ and assume that $(0:u)=0$.  Let
  $A_i=((s_iU:_Ur_i):u)$, $\forall i\in I'$. Thus $N_i=A_iF,\ \forall i\in
  I'$. By Proposition~\ref{P:compact}
  $\mathcal{F}'$ has a non-empty intersection. So $\mathcal{F}$ has
  a non-empty intersection too. 

$(2)$. We apply $(1)$ by taking $U=R_J$. \qed

\bigskip
\textbf{Proof of Theorem~\ref{T:main}}

Let $J$ be a prime ideal and $E$ a module. If $E$ is fp-injective, $E$ is a pure submodule of an injective module $M$. It follows that
$E_J$ is a pure submodule of $M_J$. So, if $M_J$ is injective we conclude that $E_J$ is fp-injective. Now we assume that $E$ is injective.

$(1)$. Suppose that $J\subset Z$. Let $s\in Z\setminus J$. Then there exists $0\ne r\in J$ such that $sr=0$. Hence $rE$ is contained
in the kernel of the natural map: $E\rightarrow E_J$. Moreover $R_J=(R/rR)_J$ and $E_J=(E/rE)_J$. By Proposition~\ref{P:prinj}, $E/rE$ is injective 
over $R/rR$ and by \cite[Theorem 11]{Cou} $R/rR$ is an IF-ring. So $E/rE$ is flat over $R/rR$. From Proposition~\ref{P:tenspur} we deduce that
$E_J$ is pure-injective and by \cite[Proposition 1.2]{Cou1} $E_J$ is fp-injective. So $E_J$ is injective.

Assume that $Z\subset J$. We set 
\[F=\{x\in E\mid J\subset (0:x)\}\qquad \mathrm{and}\qquad G=\{x\in E\mid J\subseteq (0:x)\}.\] 
Let $x\in E$ and $s\in R\setminus J$
such that $sx\in F$ (respectively $G$). Then $sJ\subset (0:x)$ (respectively $sJ\subseteq (0:x)$). Since $s\notin J$ we have $sJ=J$. Consequently
$x\in F$ (respectively $G$). Thus the multiplication by $s$ in $E/F$ (and $E/G$) is bijective because $E$ is injective. So $E/F$ and $E/G$ are 
modules over $R_J$ and $E_J\cong E/F$. We have $G\cong\mathrm{Hom}_R(R/J,E)$. It follows that $E/G\cong\mathrm{Hom}_R(J,E)$. But $J$ is a flat
module. Thus $E/G$ is injective. Let $A$ be an ideal of $R_J$ and $f:A\rightarrow E/F$
 an homomorphism. Then there exists an homomorphism $g:R_J\rightarrow E/G$ such that $g\circ u=p\circ f$ where $u:A\rightarrow R_J$ and 
$p:E/F\rightarrow E/G$ are the natural maps. It follows that there exists an homomorphism $h:R_J\rightarrow E/F$ such that $g=p\circ h$. It is easy to check
that $p\circ (f-h\circ u)=0$. So there exists an homomorphism $\ell:A\rightarrow G/F$ such that $v\circ\ell=f-h\circ u$ where $v:G/F\rightarrow E/F$ is 
the inclusion map. First assume that $A$ is finitely generated over $R_J$. We have $A=R_Ja$. If $0\ne\ell(a)=y+F$, where $y\in G$, then $(0:a)\subseteq Z\subseteq J=(0:y)$.
Since $E$ is injective there exists $x\in E$ such that $y=ax$. Hence $f(a)=a(h(1)+(x+F))$. Now suppose that $A$ is not finitely generated over $R_J$. If $a\in A$ then
there exist $b\in A$ and $r\in J$ such that $a=rb$. We get that $\ell(a)=r\ell(b)=0$. Hence $f=h\circ u$.

$(2)$. Let the notations be as above. Then $E_Z=E/F$. If $Z$ is flat, we do as above to show that $E_Z$ is injective. If $E$ is flat then $F=0$, whence $E_Z=E$.
Now, assume that $E_Z$ is fp-injective and  $Z$ is not flat. By \cite[Theorem 10]{Cou} $R_Z$ is an IF-ring. It follows that $E_Z$ is flat. Consequently $F$ is a pure 
submodule of $E$.  Suppose there exists $0\ne x\in F$. If $0\ne s\in Z$ then $(0:s)\subseteq Z\subset (0:x)$. So, there exists $y\in E$ such that $x=sy$. By \cite[Lemma 2]{Cou} 
 $(0:y)=s(0:x)\subseteq Z$. Since $F$ 
is a pure submodule, we may assume that $y\in F$. Whence $Z\subset (0:y)$. We get a contradiction. Hence $F=0$ and $E$ is flat. \qed

\bigskip Now we give a consequence of Theorem~\ref{T:main}. Recall that a 
 domain $R$ is said to be \textbf{h-local} if $R/I$ is semilocal for every
nonzero ideal $I,$ and if $R/P$
is local for every nonzero prime ideal $P,$ \cite{Mat}.

\begin{corollary} \label{C:h-loc}
Let $R$ be a h-local Pr\"ufer domain. For each multiplicative subset $S$ of $R$ and for every injective $R$-module $E$, $S^{-1}E$ is injective.
\end{corollary}
\textbf{Proof.} By \cite[Theorem 24]{Mat} $E_P$ is injective for each maximal ideal $P$. Since $R_P$ is a valuation domain, we deduce from Theorem~\ref{T:main}
 that $E_J$ is injective for each prime ideal $J$. It is easy to check that $S^{-1}R$ is a h-local Pr\"ufer domain. So, by \cite[Theorem 24]{Mat} $S^{-1}E$ is injective. \qed

\end{document}